\documentstyle[12pt,amstex,leqno]{amsart}

\def\ds{\displaystyle}
\newcommand{\R}{{\Bbb R}}
\def\tl{\widetilde}
\def\ol{\overline}
\def\sm{\setminus}
\def\B{{\cal B}}
\def\I{{\cal I}}
\def\P{{\cal P}}
\def\L{{\cal L}}
\def\NN{{\cal N}}
\def\V{{\cal V}}
\def\W{{\cal W}}
\def\l{\ell}

\newtheorem{thm}{Theorem}[subsection]  
\newtheorem{lem}[thm]{Lemma}	       
\newtheorem{crlr}[thm]{Corollary}      
\newtheorem{prp}[thm]{Proposition}     

\makeatletter
\def\theequation{\thesection.\@arabic\c@equation}
\def\thethm{\@arabic\c@thm}

\begin{document}

\title[$W^{2,1}_p$ SOLVABILITY FOR PARABOLIC POINCAR\'E PROBLEM
]{$ W^{2,1}_p$ SOLVABILITY FOR PARABOLIC POINCAR\'E PROBLEM }

\author [L.G. Softova]{L.G. Softova}
\address{ Bulgarian Academy of Sciences,
Institute of Mathematics and Informatics,
Acad. G. Bonchev Str. bl. 8, 1113 Sofia, Bulgaria}
\curraddr{Department of Mathematics, University of Bari, 4~E. Orabona~Str.,
 70~125~Bari,	Italy}
\email{luba{@}dm.uniba.it}

\renewcommand{\baselinestretch}{1.5}

\maketitle

\section{Introduction}\label{s1}
\setcounter{equation}{0}

Let $\Omega\subset \R^n,$ $n\geq 3,$ be  a bounded domain and $\L$ be an
elliptic second order differential operator defined in $\Omega.$ Consider
a vector field $\l(x)$	defined on $\partial \Omega.$ It generates a first
order boundary  operator $\B$  throught oblique derivative with respect to $\l(x).$
The couple  $(\L,\B)$  defines an elliptic {\it oblique derivative
problem (ODP).\/}  It is well  known that a boundary value problem for any 
elliptic (and parabolic) differential operator is well-posed if it satisfies the 
{\it Shapiro-Lopatinskii (S-L) complementary condition.\/} In a case of a second order 
operator  the ODP is {\it regular} (or it satisfies the mentioned  condition) when the generating  vector field  is nowhere tangential to the boundary. Otherwise the problem 
is called {\it degenerate\/} or {\it tangential.\/} An exception of this rule is the two dimensional case when the
 {\it S-L condition} always holds, i.e. the oblique derivative problem in the plane
 is always regular even when $\l(x)$ is tangential in some points of the boundary.
 The theory of regular	elliptic and parabolic  ODP  elaborated in
 H\"older and Sobolev spaces could  be found in
\cite{ADN}, \cite{H}, \cite{L}, \cite{MPS} and many others.

 For a first time a tangential problem for the Laplace
operator $(\Delta,\B)$ was posed   by Poincar\'e   in
his study of the theory of tides, but he did not solve it. Later it became an object
 of study through   modern mathematical techniques as  the theory of
pseudo differential operators, the Fourier integral operators and others.
The Poincar\'e problem still provocates interest becous of   different effects
that occure near to the set of tangency and the lack of unify approach in its study.

The qualitative properties of the solution strongly  depend on the
behavior  of the $\l$-integral curves near to the set of tangency  (order of contact,   
direction of integral curves, etc.).   In this
 context we can distinguish three	  types of contact:
 {\it neutral\/} --- the $\l$-curves  always enter into  or
leave  $\Omega$ (i.e. $l$
 preserves its sign near to the set of tangency); {\it emergent\/} --- the $\l$-curves first enter
into  $\Omega$ and after the  contact  go out  of
 the domain
($\l$ changes its sign from $-$ to $+$ near to the set of tangency); {\it
submergent\/} --- al contrary of the previous case $\l$
changes its sign from $+$ to $-$  (under a sign of $\l$ we always
mean  the sign of the scalar product $(\l\cdot\nu)$ and $\nu$ is the unit outer normal to $\partial\Omega$).  A various results about
elliptic and
 parabolic Poincar\'e  problems in H\"older and $H^s$ spaces   are
 presented in \cite{Mz}, \cite{P}, \cite{Pl} \cite{PP}, \cite{W1}, \cite{W2}, \cite{W3},
\cite{MPV}.

In the present work  we consider degenerate  ODP in a bounded cylinder 
 $Q=\Omega\times(0,T)$ for a linear second-order
 uniformly parabolic
operator $\P$ with coefficients allowing discontinuity in $t.$ 
The vector field
$\l(x,t)$ generating   $\B$ is	 defined on  $S=\partial\Omega\times(0,T)$
 and is tangential to it in some subset $E.$ The kind of contact
is of neutral type and we suppose that $\gamma(x,t)=(\l(x,t)\cdot\nu(x))\geq 0$ on $ S.$
 It means that	the boundary value problem under
 consideration is of Fredholm type,
 i.e.  both the kernel and cokernel  are of finite dimension.

We are interested of strong solvability of our problem in $W^{2,1}_p(Q),$ $p\in (1,\infty).$
Because of the loss of regularity of the solution near to the set of tangency $E$ we impose higher regularity in $E$ of the data.       
The study is based on the original Winzel's idea to extend $\l$ into $\Omega$ such
 that to obtain explicit representation of the solution through the integral curves
of that
 extension. Thus the problem is reduced  to obtaining of suitable a priori estimates for the solution and its derivatives on an expanding family of cylinders.  
Further, the solvability is proved using    regularization technique
 which,  roughly speaking, means to perturb the vector field $\l$ by adding
small $\varepsilon$ times $\nu,$ to solve the such obtained {\it regular\/} ODP
  and then pass to limit as $\varepsilon\to 0.$ The perturbed problem
regards linear  uniformly parabolic
  operator $\P$ with VMO coefficients  and boundary
 operator
	$\B$
with $(\l_\varepsilon\cdot\nu)>0.$	In this case we dispose of  unique  solvability result in
 $W^{2,1}_p(Q),$ $p\in(1,\infty)$ supposing $\P u\in L^p(Q)$
 and  initial and  boundary data  belonging to the corresponding Besov
 spaces
(see \cite{Sf}, \cite{MPS}).

Poincar\'e problem for linear uniformly parabolic operators with H\"older
 continuous coefficients is studied in \cite{Pl} (see also \cite{PP})
 where unique solvability in
the corresponding H\"older spaces is obtained.  Moreover, the linear results
were applied to the study of
semilinear parabolic problem in H\"older spaces.
A tangential ODP for second-order	uniformly elliptic operators
with Lipschitz
  continuous coefficients was studied in \cite{MPV} (see also \cite{MPS}). It is
obtained
 strong solvability  in $W^{2,p}(\Omega)$ but for $p>n/2.$

 In our case
the parabolic
 structure of the  equation permits  to obtain an a priori
 estimate for the solution in $W^{2,1}_p(Q)$  only through the data of the
problem. Thus we are able to  prove
 unique solvability for all $p\in(1,\infty)$ avoiding the use of maximum
principle and omitting any additional conditions on the vector field.

\section{Statement of the problem and main results}\label{s2}
 \setcounter{equation}{0}
\setcounter{thm}{0}

Let $\Omega\subset \R^n,$ $n\geq 3$ be a bounded domain with
$\partial\Omega \in C^{2,1}
$  and	$Q=\Omega\times(0,T)$ be	 a cylinder in $\R^{n+1}.$  Set
$\l(x,t)=(\l^1(x,t),\ldots,\l^n(x,t),0)$ for a unit  vector field  defined on the
lateral boundary	    $S=\partial\Omega\times(0,T).$
We consider the following  {\it oblique derivative problem\/}
 \begin{equation}\label{1}
\begin{cases}
\P u\equiv u_t - a^{ij}(x,t)D_{ij}u=f(x,t)& {\text in }\  Q,\\
\I u \equiv u(x,0)=\psi(x) &	{\text on}\	\Omega,\\
\B u \equiv\ds \frac{\partial u}{\partial\l}= \l^i(x,t)D_iu=\varphi (x,t) &
{\text on}\  S.
\end{cases}
\end{equation}
Denote by $\nu(x)=(\nu^1(x),\ldots,\nu^n(x))$  the unit outward normal to
$\partial\Omega.$ Then we can write $\l(x,t)=\tau(x,t)+\gamma(x,t)\nu(x)$ where
$\tau(x,t)$ is	tangential projection of $\l$ on $S$ and
$\gamma(x,t)=(\l(x,t)\cdot\nu(x))\geq 0.$
Let $E\subset S$ be the set of tangency and $E\cap \partial\Omega=E_0.$

The set of tangency has the form $E=E_0\times(0,T)$ where $E_0\subset \partial\Omega.$
By $\Sigma$ we denote a cylinder with a base $\Sigma_0\subset\Omega$ being small
neighborhood 
of $E_0,$  $S_\Sigma=\partial \Sigma\cap S$ such that
$E\subset S_\Sigma.$ Thus $\Sigma=\Sigma_0\times(0,T)$ is a domain where we shall impose
more restrictive conditions on the data of \eqref{1},   while
	in   $Q\sm\Sigma$ we can take  the
same   conditions as in the regular  case.  Precisely:
\begin{enumerate}
\item[$(i)$]  $\P$ is a {\it uniformly parabolic operator\/}:
$\exists$  $\lambda>0$ such that
$$
\begin{cases}
\lambda |\xi|^2\leq a^{ij}(x,t)\xi_i\xi_j\leq \lambda^{-1}|\xi|^2 &
 {\rm a.a.}\
(x,t)\in Q,\  \forall \xi\in\R^n,\\
a^{ij}\in V^1_\infty(\Sigma)\cap  VMO(Q), & a^{ij}=a^{ji}\quad(\Longrightarrow
 a^{ij}\in L^\infty(Q)).
  \end{cases}
$$
Here $VMO$ is the Sarason class of function with {\it vanishing mean oscillation}
(see \cite{S}) and $V^1_\infty (\Sigma):=\{v, D_x v\in L^\infty(\Sigma)\}.$
\item[$(ii)$]  $\B$ is a {\it degenerate oblique derivative operator\/} defined
through a tangential vector field of {\it neutral\/} type:
$$
\begin{cases}
 &\gamma(x,t)=(\l(x,t)\cdot \nu(x))\geq 0  \ {\rm on} \  S,\quad
\gamma(x,t)=0\ {\rm on}\
 E\subset S,\\
&\l^i(x,t) \in	{\rm Lip}(S)\cap W^{2,1}_\infty(S_\Sigma)
\end{cases}
$$
where the space $W^{2,1}_\infty(S_\Sigma)$  consists of functions having
 $L^\infty(S_\Sigma)$ derivatives in $x$ up to order $2$ and in $t$ of order $1.$    
\item[$(iii)$] {\it Regularity of the data}: for all $p\in (1,\infty)$
\begin{enumerate}
\item[$(iii_a)$]
$
f(x,t)\in \V_p(Q):=L^p(Q)\cap	V^1_p(\Sigma),$\\
$\|f\|_{\V_p(Q)}=\|f\|_{p,Q}+\|f\|_{V^1_p(\Sigma)}=\|f\|_{p,Q}+\|D_xf\|_{p,\Sigma};$
\item[$(iii_b)$]
$\psi(x)\in W^{3-2/p}_{p}(\Omega),$ 
$
\ds\|\psi\|^{(3-2/p)}_{p,\Omega}=\sum_{s=0}^{[3-2/p]}\sum_{(s)}\|D^s_x\psi
\|_{p, \Omega } + \ll \ds\psi\gg^{(3-2/p)}_{p,\Omega}$\\
$\ds \ll \psi\gg^{(3-2/p)}_{p,\Omega}=\sum_{(s=[3-2/p])}\left(\int_\Omega
dx\int_\Omega |D^s_x\psi(x)-D^s_y\psi(y)|^p \frac{dy}{|x-y|^{n+p\alpha}}
\right)^{1/p}$ where $\alpha=3-2/p- [3- 2/p].$\\[5pt]
\item[$(iii_c)$]
$ \varphi(x,t) \in  \W_p(S):= W^{1-1/p,1/2-1/2p}_p(S)\cap
 W^{2-1/p,1-1/2p}_p(S_\Sigma),$\\
$ \|\varphi\|_{\W_p(S)}= \|\varphi\|_{p,S}^{(1-1/p)}
+\|\varphi\|_{p,S_\Sigma}^{(2-1/p)},$ where
$$
\|\varphi\|^{(l)}_{p,S}=\sum_{0\leq 2r+s\leq [l]}\|D^r_tD^s_x \varphi
\|_{p,S} + \ll \varphi\gg^{(l)}_{p,S}
$$
$$
\ll \varphi\gg^{(l)}_{p,S}
=\sum_{2r+s=[l]}\ll D_t^rD^s_x \varphi
\gg^{(l-[l])}_{p,x;S}
+ \sum_{0<l-2r-s<2} \ll D_t^rD^s_x \varphi \gg^{((l-2r-s)/2)}_{p,t;S},
$$
and for $0<\alpha<1$ we define
$$
\ll v\gg^{(\alpha)}_{p,x;S}=\left(\int_0^T
dt\int_{\partial\Omega}
dx \int_{\partial\Omega} |v(x,t) -v(y,t)|^p
\frac{dy}{|x-y|^{n-1+p\alpha}} \right)^{1/p},
$$
$$
\ll \varphi\gg^{(\alpha)}_{p,t;S}=\left(\int_{\partial\Omega}
dx\int_0^T dt  \int_0^T   |v(x,t) -v(x,\tau)|^p
\frac{d\tau}{|t-\tau|^{1+p\alpha}} \right)^{1/p}.
$$
\item[$(iii_d)$]  {\it Compatibility condition } on $\partial\Omega:$
$$
\B \psi(x)=\varphi(x,0)\quad {\rm for}\quad \begin{cases}
 x\in E_0 & {\rm and}\quad p>3/2,\\
 x\in \partial\Omega\sm E_0 & {\rm and}\quad p>3.
\end{cases}
$$
\end{enumerate}
\item[$(iv)$] The integral curves of $\l$ on $E$ are {\it  non closed\/}
and of {\it finite length.\/}
 \end{enumerate}

We are interested in  solvability of the problem  \eqref{1} in
the Sobolev space
$$
W^{2,1}_p(Q)=\left\{u\in L^p(0,T;W^{2,p}(\Omega)),u_t\in L^p(Q),p\in(1,\infty) \right\}
$$
endowed by the norm
\begin{equation}\label{Sobolev}
\|u\|_{W^{2,1}_p(Q)}=\sum_{j=0}^2 \ll u\gg^{(j)}_{p,Q}	=\sum_{j=0}^2 \sum_{(2r+s=j)}\|D^r_tD^s_xu\|_{p,Q}.
\end{equation}
Under a {\it strong solution\/} to \eqref{1} we mean a function $u\in W^{2,1}_p(Q)$
satisfying $\P u=f$ almost everywhere in $Q$ and  the initial and
boundary conditions  hold in trace sense.

\begin{thm}[A priori estimate]\label{th2}
   Suppose conditions $(i)-(iv)$  to be
 fulfilled, and $u\in  W^{2,1}_p(Q)$  for $1<p<\infty$.  Let  $\P u\in
\V_p(Q),$ $\I u\in W^{3-2/p}_p(\Omega)$  and $\B u\in
\W_p(S)$  then
\begin{equation}\label{eq5}
\|u\|_{W^{2,1}_p(Q)} \leq
C\left( \|\P u\|_{\V_p(Q)}+\|\I u\|_{W^{3-2/p}_p(\Omega)}+\|\B u\|_{\W_p(S)}\right)
\end{equation}
where the constant depends on
$n,p,\lambda,T,\partial\Omega,\l, \|a^{ij}\|_{L^\infty(\Sigma)},$ and
 the $VMO$-moduli of the coefficients.
\end{thm}
\begin{thm}[Unique strong solvability]\label{th3}
   Assume $(i)-(iv) $to be
fulfilled.
Then for all $f\in \V_p(Q),$ $\psi\in W^{3-2/p}_p(\Omega)$ and
$\varphi\in \W_p(S),$ the problem
(\ref{1}) admits a unique  solution $u\in W^{2,1}_p(Q)$ for all
$p\in(1,\infty).$
\end{thm}

As in the case of regular ODP,	the  embedding result
\cite{LSU} gives   H\"older continuity of the solution
to (\ref{1}) for appropriate values of $p.$
\begin{crlr}\label{crlr1}
Let $u\in W^{2,1}_p(Q)$ be a solution of \eqref{1}. Then
\begin{itemize}
\item[$i)$] $u\in C^{0,\alpha}(\bar Q)$ with $\alpha<2-(n+2)/p$ if
$p\in((n+2)/2, n+2];$
\item[$ii)$] $D_xu\in C^{0,\beta}(\bar Q)$ with $\beta<1-(n+2)/p$ if $p>n+2.$
\end{itemize}
\end{crlr}

\section{Auxiliary results}\label{s5}

The following assertion gives some
 geometrical properties of  $\l$  (see \cite{W3}, \cite{PP}).
\begin{prp}\label{Geometry}
 Let $\l$ and $\Omega$ satisfy the assumptions listed above.
 There exists a finite upper bound $\kappa_0$ of the arclength of
$\l$-integral curves lying on the set $E.$ Moreover,
there
exist extensions $L(x,t)\in  W^{2,1}_\infty(\Sigma)$ and
 $\bar\nu(x)\in C^{1,1}(\Sigma_0)$ of $\l(x,t)$ and  $\nu(x),$ and a cylinder $Q_0=\Omega_0\times(0,T)$
 with the following properties:

The base $\Omega_0\subset\Omega,$ $\partial\Omega_0\subset C^{1,1}$ is such that 
$E_0\subset \partial\Omega\sm \partial\Omega_0.$ Denote by $\partial Q_0$ the lateral
 boundary of $Q_0$ and $S_0:=\partial Q_0\sm S.$ The extension $L(x,t)$ is strictly 
transverzal to $S_0$ and each point of $S_\Sigma$ can be reached from $S_0$ along  an
$L$-integral curve of length at most $\kappa'>\kappa_0.$

Define $
Q_\tau=Q_0\cup \{ e^{sL}(x,t)\in Q:\
(x,t)\in S_0, \  0\leq s\leq \tau,\ \tau>0\}.
$
Under  $\partial Q_\tau$ we understand  the lateral boundary of  $Q_\tau,$  and  $S_\tau=\partial
 Q_\tau\sm S=\{{e^{\tau L}(x,t),(x,t)\in S_0}\}$
The  family $\{Q_\tau\}_{\tau\geq 0}$ is non-decreasing  and
 for every $\delta>0$ there exists $\theta=\theta(\delta)>0,$
independent of $\tau,$	such that	${\rm dist\,}(S_\tau,
S_{\tau+\delta})\geq \theta$ whenever $Q\setminus
Q_{\tau+\delta}\not= \emptyset.$
 The field $L$ is strictly transversal to
$S_\tau  \in W^{2,1}_\infty$ uniformly in  $\tau.$
\end{prp}
Let we note that in the geometrical construction above we want that
 $E_0\subset \partial \Omega\sm\partial\Omega_0$ although the set 
$\Omega_0\cap\Sigma_0$ could not be empty. In fact in the cylinder 
$Q_0=\Omega\times(0,T)$ we have regular ODP with boundary vector field 
$L$ which coincides with $\l $ on $\partial Q_0\cap S.$

It is well known (see \cite{GT}) that there exists a
neighborhood $\NN$ of $\partial\Omega$ such that  for any $x\in\NN$
 there exists unique closest point $y(x)\in\partial\Omega$
and  $x=y(x)-\nu(y)d(x),$  $y,d\in C^{1,1}(\bar \Omega).$
Let $\ol{\Omega\sm\Omega_0}\subset\NN,$ then $d,y\in C^{1,1}
(\ol{\Omega\sm\Omega_0})$
 and we set
$$
 \bar\nu(x)=\nu(y(x)),\quad L(x,t)=\l(y(x),t)+d(x)\bar\nu(x)\quad \forall x\in
\ol{\Omega\setminus\Omega_0}, t\in(0,T).
$$
The regularity of $L$ follows by the regularity
properties of $\l$ and	$d.$
The rest of the proof repeats the arguments in
 \cite[Proposition~3.2]{W3} and \cite[Proposition~3.2.5]{PP}.

We  need also of the following	variant of the	Gronwall
inequality (see \cite{W3}).
\begin{prp}\label{prop2}
Let $\zeta(\tau)$ be continuous, bounded and positive function, defined on
$[0,\infty).$  Suppose there exist positive constants $\delta,$ $A$ and $C$ such
that
$$
\zeta(\tau)\leq A+C\int_0^\tau \zeta(s+\delta) ds\quad{\rm for\  all}\quad \tau>0.
$$
Then   
$
\zeta(\tau)\leq A\left(1+(2\pi)^{-1/2}(1-C\delta e)^{-1}e^{\tau/\delta} \right)
$ for $C\delta e<1.$
\end{prp}

\section{Cauchy problem}\label{s6}

In the present section we discuss the solvability of the Cauchy problem
and the  reducing of \eqref{1} to one with homogeneous initial data. Consider
\begin{equation}\label{CP}
\begin{cases}
v_t-\Delta v=0 & (x,t)\in \R^n\times (0,T)\\
v(x,0)=\tl\psi(x) &  x\in \R^n
\end{cases}
\end{equation}
where $\tl\psi\in W^{3-2/p}_p(\R^n)$ is an extension of $\psi$  as
zero for $x\not\in\bar\Omega.$	 The solution of \eqref{CP} is given by the potential
$$
v(x,t)=(\Gamma\ast_1\tl\psi)=\frac{1}{(4\pi t)^{n/2}}\int_{\R^n}
e^{-\frac{|x-y|^2}{4t}}\tl\psi(y)dy =\frac{1}{(4\pi t)^{n/2}}\int_{\Omega}
e^{-\frac{|x-y|^2}{4t}}\psi(y)dy
$$
and it can be considered as an
 extension in $t>0$
of the initial data  preserving its
regularity.
In our case we have supposed higher regularity of $\I u$ (note that when
 the boundary  operator is regular  it is enough to take
$\I u \in W^{2-2/p}_p(\Omega)$) looking for a solution of \eqref{CP}
 possessing higher
regularity in $x.$
Having in mind the estimates of the heat potential in Sobolev spaces
(see \cite[Ch~IV, $\S~3,4$]{LSU})  it follows that $v\in
W^{2,1}_p(\R^n\times(0,T))$
 and   corresponding a priori estimate holds.
 Moreover, since
$D_xv=D_x(\Gamma\ast_1\tl\psi)=(\Gamma\ast_1 D_x\tl\psi)$ we obtain also that
 $$
 \ll(\Gamma\ast_1 D_x\tl\psi)\gg^{(2)}_{p,\R^n\times\R_+}\leq C
\ll D_x\tl\psi\gg^{(2-2/p)}_{p,\R^n}= \ll\tl\psi\gg^{(3-2/p)}_{p,\R^n}.
$$
On the other hand according to \eqref{Sobolev} and $(iii_b)$ we can write
\begin{align*}
&\ll(\Gamma\ast_1 D_x\tl\psi)\gg^{(2)}_{p,\R^n\times\R_+} =
\|D^2_x(\Gamma\ast_1 D_x\tl\psi)\|_{p,\R^n\times\R_+} +
\|D_t(\Gamma\ast_1 D_x\tl\psi)\|_{p,\R^n\times\R_+}\\
&\quad=\|D^3_x v\|_{p,\R^n\times\R_+} +
\|D_tD_x v\|_{p,\R^n\times\R_+}\leq C \ll\tl\psi\gg^{(3-2/p)}_{p,\R^n}
\leq C\|\psi\|^{(3-2/p)}_{p,\Omega}.
\end{align*}

The function $w=u-v$ satisfies
 $$
\begin{cases}
\P w\equiv w _t - a^{ij}D_{ij}w=f(x,t) +(a^{ij}-\delta^{ij})D_{ij}v :=\tl f(x,t)&
{\text in }\  Q,\\
\I w \equiv w(x,0)=0 &	  {\text on}\	  \Omega,\\
\B w \equiv\ds \frac{\partial w}{\partial\l}= \varphi
(x,t)+\l^i(x,t)D_iv:=\tl\varphi(x,t) & {\text on}\  S
\end{cases}
$$
where $\tl f\in \V_p(Q)$ and $\tl\varphi\in \W_p(S)$ according to $(iii)$  and	the
regularity  of  $v.$

\section{A priori estimates}\label{s3}

In our further considerations we always  suppose that $\I u=0.$
 Our main goal is to obtain an a
priori estimate for the solution  only through the data  that allows  proving of  unique solvability of \eqref{1}.
\begin{lem}\label{th1}
Let $v\in W^{2,1}_p(Q),$ $p\in(1,\infty),$ $v(x,0)=0$ and $(i)$ hold.	 Let
$\Omega_1\subset\Omega_2\subset\Omega$	with $\partial\Omega_i\in C^{1,1},$
 $Q_i=\Omega_i\times(0,T),$ $S_i=\partial Q_i\cap  S,$
 $i=1,2$ such that ${\rm
dist\,}(Q_1,Q\sm Q_2)\geq \theta>0.$ Then there
exist constants $C'$ and $C''(\theta)$ such that
$$
\|v\|_{W^{2,1}_p(Q_1)}\leq  C'\left(\|\P v\|_{p,Q_2} +
\|v\|^{(2-1/p)}_{p,S_2} \right)
+C''(\theta)\left(\|v\|_{V^1_p(Q_2)} +
\|v\|^{(1-1/p)}_{p,S_2}\right).
$$
\end{lem}
\begin{pf}
Take a cutoff function	$\eta(x)\in C^\infty_0(\R^{n})$
 such that  $\eta(x)=1$ for $ x\in\bar\Omega_1,$ $\eta(x)=0$ for $x\in
\ol{\Omega\sm\Omega_2}$
and $\max_{\bar\Omega}|D^\alpha\eta(x)|\leq C_\alpha\theta^{-|\alpha|}.$ For
$v\in W^{2,1}_p(Q)$ it is easy to see that $\P (v\eta)$ belongs to $L^p(Q).$
Having in mind	the
a priori estimate
  for homogeneous Cauchy-Dirichlet problem  (see \cite{BC}),  it is a standard
procedure    to obtain an analogous estimate for non
 homogeneous one
\begin{align*}
\|v\|_{W^{2,1}_p(Q_1)} & \leq	 \|\eta v\|_{W^{2,1}_p(Q_2)}
\leq  \left(\|\P (\eta v)\|_{p,Q_2} +\|\eta v\|^{(2-1/p)}_{p,\partial Q_2}\right)\\
& \leq
C_0 \left(\|\P v\|_{p,Q_2} +
\|v\|^{(2-1/p)}_{p,S_2} \right) + \frac{C_1}{\theta} \left(\|Dv\|_{p,Q_2}  +
\|v\|^{(1-1/p)}_{p,S_2} \right) + \frac{C_2}{\theta^2} \|v\|_{p,Q_2} \\
& \leq	C' \left(\|\P v\|_{p,Q_2}
+\|v\|^{(2-1/p)}_{p,S_2} \right) +C''(\theta) \left(\|v\|_{V^1_p(Q_2)}
+\|v\|^{(1-1/p)}_{p,S_2} \right).
\end{align*}
\end{pf}

\noindent
{\bf Proof of Theorem~1.}
The  derivative  $\partial u/\partial L$ satisfies
in $\Sigma$ the following Cauchy-Dirichlet problem
\begin{equation}\label{CDderiv}
\begin{cases}
 \ds  D_t\big(\frac{\partial
\ds u}{\partial L}\big)-  a^{ij}D_{ij}\big(\frac{\partial u}{\partial
\ds  L}\big)= \frac{\partial f}{\partial L}+\frac{\partial
 \ds a^{ij}}{\partial L}D_{ij}u +D_ku D_tL^k&\\[5pt]
 \ds\qquad 	- a^{ij}\left(2D_{ik}u D_jL^k +D_k u
 D_{ij}L^k\right) =F(x,t) &
 {\text in }\  \Sigma\\[5pt]
 \ds   \frac{\partial u(x,0)}{\partial L}=0 & {\text on }\ \Sigma_0\\[7pt]
\ds  \frac{\partial u(x,t)}{\partial
L}=\varphi(x,t) &  {\text on }\ S_\Sigma.
\end{cases}
\end{equation}
Let $\Sigma_0''\subset\Sigma_0$  and  consider a cut-off function
 $\eta(x)\in C^2_0(\R^n)$ such that
 $\eta(x)=1$ for $x\in \bar\Sigma''_0,$  $\eta(x)=0$ for
$x\in \ol{\Omega\sm\Sigma_0}$ and $\sup_{\Sigma_0}|D^\alpha\eta(x)|\leq
Cd_0^{-|\alpha|}$
where $d_0$ is the  distance between $\Sigma_0''$ and
 $\partial \Sigma_0\sm \partial \Omega.$
Then  $V=\frac{\partial u}{\partial L}\eta$
is a solution of
$$
\begin{cases}
\ds  D_t V- a^{ij}(x,t)D_{ij}V= F(x,t)\eta(x)&\\
\qquad\  \ds- a^{ij}(x,t)\big(2D_i\eta	D_j\big(\frac{\partial u}{\partial L}\big)
\ds +  \frac{\partial u}{\partial L} D_{ij}\eta\big)=:F_1(x,t) & {\text in }\  \Sigma\\
\ds V(x,0)=0 & {\text on }\ \Sigma_0\\
\ds V(x,t) =\varphi(x,t) & {\text on }\ S_\Sigma.
\end{cases}
$$
As in Lemma~\ref{th1}  we have an a priori estimate for the solution of
  the above problem
$$
\|V\|_{W^{2,1}_p(\Sigma)}\leq
C\left(\|F_1\|_{p,\Sigma}+\|\varphi\|_{p,S_\Sigma}^{(2-1/p)} \right)\leq
  C \big(\|u \|_{W^{2,1}_p(\Sigma )}+\|f\|_{V^1_p(\Sigma)}
+\|\varphi\|_{p,S_\Sigma}^{(2-1/p)}\big)
$$
and the constant depends on the $^\infty$-norms of $a^{ij},$ $\l$ and
diam\,$\Sigma$ (through the extension $L$).
Having in mind that $V = \partial u/\partial L$ in $\Sigma''$ we obtain
\begin{equation}\label{eq2}
\big\|\frac{\partial u}{\partial L}\big\|_{W^{2,1}_p(\Sigma'')}\leq
C\left(\|f\|_{\V_p(Q)}+\|\varphi\|_{\W_p(Q)}+ \|u\|_{W^{2,1}_p(Q)}\right).
\end{equation}

Considering in analogous way the regular $ODP$
in the cylinder $Q\sm \Sigma''$ and applying   \cite[Theorem~1]{Sf} we obtain
\begin{equation}\label{eq3}
\|u\|_{W^{2,1}_p(Q\sm\Sigma'')} \leq C\left(\|f\|_{\V_p(Q)}+\|\varphi\|_{\W_p(S)}
+\|u\|_{V^1_p(Q)}\right).
\end{equation}
To estimate the Sobolev
 norm of $u$ in $\Sigma''$  we make use of 
an explicit formula for the solution through the $L$-integral curves.
Construct a cylinder  $Q_0=\Omega_0\times(0,T)$ (see	Proposition~\ref{Geometry})
 such that    $\Omega_0\subset
\{\Omega\sm\Sigma_0''\},$   $S_0$ lies in $\Sigma\sm
\Sigma''$  where  $L$ is well  defined. Let $(x,t)\in\Sigma''$ and
$\psi(\tau;x,t)=e^{\tau L}(x,t),$ $\tau\in[0,\kappa']$ be the parametrisation
of the $L$-integral curve passing through that point ($\psi(0;x,t)=(x,t)$).
	Then there
exists unique $\xi(x,t)\in W^{2,1}_\infty(\Sigma)$ such that
$\psi(-\xi(x,t);x,t)\in Q_0\subset Q\sm \Sigma''$ and
$$
u(x,t)=u\circ \psi(-\xi(x,t);x,t) +\int_{\tau-\xi(x,t)}^\tau
\frac{\partial u}{\partial L}\circ \psi(s-\tau;x,t)ds.
 $$
First we shall estimate the $L^p$-norm of $D^2u$ in $Q_\tau\cap\Sigma''$ where
$\{Q_\tau\}_{\tau\geq 0}$ is the expanding family of cylinders defined in
Proposition~\ref{Geometry} (note $(x,t)\in S_\tau$)
\begin{align}
\nonumber
\|D^2u\|_{p,Q_\tau\cap \Sigma''}&\leq \|D^2u\|_{p,Q_\tau}\leq
C\big(
\|D^2u\|_{p,Q\sm\Sigma''}+\|Du\|_{p,Q\sm\Sigma''}\\
\nonumber
+&\int_0^\tau \left(\big\|D^2\big(\frac{\partial u}{\partial L}\big)
\big\|_{p,Q_s\cap
\Sigma''}+\big\|D\big(\frac{\partial u}{\partial L}\big)
\big\|_{p,Q_s\cap \Sigma''}\right) ds
+ \big\|\frac{\partial u}{\partial L}\big\|_{V^1_p(\Sigma\sm\Sigma'')}\big) \\
\nonumber
\leq&C\big( \|u\|_{W^{2,1}_p(Q\sm\Sigma'')}
+   \int_0^\tau \big\|D^2\big(\frac{\partial u}{\partial
L}\big)\big\|_{p,Q_s\cap \Sigma''}ds
 + \big\|D\big(\frac{\partial u}{\partial L}\big)\big\|_{p,\Sigma''}\\
\nonumber
+&\big\|\frac{\partial u}{\partial L}\big\|_{V^1_p(\Sigma\sm\Sigma'')}\big) \\
\label{eq6}
\leq&C\big( \|u\|_{W^{2,1}_p(Q\sm\Sigma'')}
+   \int_0^\tau \big\|D^2\big(\frac{\partial u}{\partial
L}\big)\big\|_{p,Q_s\cap
\Sigma''}ds + \big\|\frac{\partial u}{\partial L}\big\|_{V^1_p(\Sigma)} \big)
\end{align}
We need  an upper bound for the norm under the integral in order
 to apply the Gronwall type  inequality.

In our considerations we distinguish two cases. The first one is when $Q\sm
Q_{s+\delta}\not=\emptyset$ and according to Proposition~\ref{Geometry},
  there exists	$\theta=\theta(\delta)>0$ such that ${\rm
dist\,}(S_s,S_{s+\delta}) \geq \theta.$   Consider the right-hand side in
\eqref{eq6}. The first term is estimated by \eqref{eq3}. To estimate the
second one we apply Lemma~\ref{th1} to the solution $\partial u/\partial L$
of \eqref{CDderiv} with  $Q_s\cap \Sigma''$
and $Q_{s+\delta}\cap \Sigma''$ instead of $Q_1$ and $Q_2$
\begin{align*}
 \big\|&D^2\big(\frac{\partial u}{\partial L}\big)\big\|_{p,Q_s\cap \Sigma''}
\leq \big\|\frac{\partial u}{\partial L}\big\|_{W^{2,1}_p(Q_s\cap \Sigma'')}
\leq C'\left( \|F\|_{p,Q_{s+\delta}\cap \Sigma''} + \|\varphi\|^{(2-1/p)}_{p,\partial Q_{s+\delta}\cap S}\right)\\
 &+ C''(\theta) \left(\big\|\frac{\partial u}{\partial L}\big\|_{V^1_p(Q_{s+\delta}\cap \Sigma'')} 
+\|\varphi\|^{(1-1/p)}_{p,\partial  Q_{s+\delta}\cap S} \right)\\
&  \leq C\left(\|u\|_{W^{2,1}_p(Q_{s+\delta}\cap\Sigma'')}
+\|f\|_{\V_p(Q)}+ \|\varphi\|_{\W_p(Q)} \right)
+C''(\theta)\big\|\frac{\partial u}{\partial L}\big\|_{V^1_p(Q_{s+\delta}\cap\Sigma'')}.
\end{align*}
 Let $\Sigma'_0\subset
\Sigma''_0$ and $\Sigma'=\Sigma'_0\times(0,T).$  Thus
\begin{align*}
\big\|\frac{\partial u}{\partial L}& \big\|_{V^1_p(Q_{s+\delta}\cap \Sigma'')}
\leq  \big\|\frac{\partial u}{\partial L}\big\|_{V^1_p(\Sigma')}
+\big\|\frac{\partial u}{\partial L}\big\|_{V^1_p(\Sigma''\sm\Sigma')}\\
 &\leq  \varepsilon \big\|D^2\big(\frac{\partial u}{\partial L}\big)\big\|_{p,\Sigma'}+
C(\varepsilon) \big\|\frac{\partial u}{\partial L}\big\|_{p,\Sigma'}+
C\|u\|_{W^{2,1}_p(Q\sm\Sigma')}\\
 &\leq  \varepsilon\big\|\frac{\partial u}{\partial L}\big\|_{W^{2,1}_p(\Sigma')}
 +C(\varepsilon) \big\|\frac{\partial u}{\partial L}\big\|_{p,\Sigma'}
+C( \|f\|_{\V_p(Q)}+\|\varphi\|_{\W_p(Q)} +\|u\|_{V^1_p(Q)}) \\
& \leq  C( \|f\|_{\V_p(Q)}+\|\varphi\|_{\W_p(Q)}) +C(\varepsilon)\|u\|_{V^1_p(Q)}
+\varepsilon\|u\|_{W^{2,1}_p(\Sigma')}
\end{align*}
after applying the Gagliardo-Nirenberg
interpolation  inequality with suitable $\varepsilon>0,$ \eqref{eq2} and  \eqref{eq3}.
Finally, the last term in \eqref{eq6} is estimated dividing it in two
parts, i.e.
\begin{align}\label{eq18}
\nonumber
\big\|\frac{\partial u}{\partial L} \big\|_{V^1_p(\Sigma'')}& +
\big\|\frac{\partial u}{\partial L}\big\|_{V^1_p(\Sigma\sm\Sigma'')}\leq
 \varepsilon \big\|D^2\big(\frac{\partial u}{\partial L}
\big)\big\|_{p,\Sigma''}\\
\nonumber
& +C(\varepsilon)\big\|\frac{\partial u}{\partial
L}\big\|_{p,\Sigma''} + \big\|\frac{\partial u}{\partial
L}\big\|_{V^1_p(\Sigma\sm\Sigma'')}\\
&\leq \varepsilon \|u\|_{W^{2,1}_p(Q)} + C\big(\|f\|_{\V_p(Q)}
+\|\varphi\|_{\W_p(Q)}\big)  +C(\varepsilon)\|u\|_{V^1_p(Q)}
\end{align}
using  \eqref{eq2} and \eqref{eq3}	in	the last step.
Substituting the
above estimates in \eqref{eq6} we obtain
\begin{align} \label{eq12}
\nonumber
 \|D^2u\|_{p,Q_\tau\cap \Sigma''}&\leq
C\left(\|f\|_{\V_p(Q)}+\|\varphi\|_{\W_p(Q)} \right)
+\varepsilon\|u\|_{W^{2,1}_p(Q)}\\
&+C \int_0^\tau\|u\|_{W^{2,1}_p(Q_{s+\delta}\cap\Sigma'')} ds
 +C(\theta,\varepsilon)\| u\|_{V^1_p(Q)}.
\end{align}
From  the equation $u_t=a^{ij}D_{ij}u+f$ it follows
an analogous  estimate also for  $u_t.$ Hence
\begin{align}\label{eq13}
\nonumber
\|D^2u\|_{p,Q_\tau\cap \Sigma''}&+\|u_t\|_{p,Q_\tau\cap \Sigma''}\leq
C\big(\|f\|_{\V_p(Q)}+\|\varphi\|_{\W_p(Q)}\big)\\
\nonumber
&+\varepsilon\|u\|_{W^{2,1}_p(Q)}+ C \int_0^\tau\big( \|D^2u\|_{p,Q_{s+\delta}\cap\Sigma''}
+\|u_t\|_{p,Q_{s+\delta}\cap\Sigma''} \big) ds\\
& +C(\theta,\varepsilon)\| u\|_{V^1_p(Q)}
\end{align}
where we have interpolated the norms of the lower derivatives of $u$ 
under the integral and use that $\int_0^\tau \|u\|_{p, Q_{s+\delta}\cap\Sigma''}ds\leq C\|u\|_{p,Q}.$

In the second case  $Q_{s+\delta}$ covers the whole cylinder and hence  $Q\sm Q_{s+\delta}=\emptyset.$  The
difference with the first one is the estimate  for the second term in
\eqref{eq6}. Using \eqref{eq2} and \eqref{eq3} instead of
Lemma~\ref{th1} and having in mind $\Sigma''=Q_{s+\delta}\cap \Sigma'',$ we get
\begin{align}\label{eq19}
\nonumber
\big\|\frac{\partial u}{\partial L} \big\|_{W^{2,1}_p(Q_s\cap\Sigma'')}
&\leq  C\big(\|f\|_{\V_p(Q)} +\|\varphi\|_{\W_p(Q)}\\
\nonumber
&+\|u\|_{W^{2,1}_p(Q_{s+\delta}\cap \Sigma'')}+ \|u\|_{W^{2,1}_p(Q\sm\Sigma'')} \big)\\
\nonumber
&\leq  C\left(\|f\|_{\V_p(Q)} +\|\varphi\|_{\W_p(Q)} + \|u\|_{V^1_p(Q)}\right)\\
&+\|u\|_{W^{2,1}_p(Q_{s+\delta}\cap \Sigma'')}.
\end{align}
Finally estimating  \eqref{eq6} through \eqref{eq3}, \eqref{eq19}  and \eqref{eq18}
we obtain again \eqref{eq12} and \eqref{eq13}.

Now we can apply the Proposition~\ref{prop2} to the function
$\zeta(\tau)=\|D^2u\|_{p,Q_\tau\cap \Sigma''}+ \|u_t\|_{p,Q_\tau\cap \Sigma''}.$
Hence for values of $\tau$ grater than $\kappa'$ for which 
 $ Q_\tau\cap \Sigma''\equiv  \Sigma''$
and $\delta>0$ small enough  we have
\begin{align}\label{eq4}
\nonumber 
\|D^2u\|_{p, \Sigma''}&+\|u_t\|_{p, \Sigma''}\leq
C(\|f\|_{\V_p(Q)}+\|\varphi\|_{\W_p(Q)})\\
&+\varepsilon\|u\|_{W^{2,1}_p(Q)} + \|u\|_{V^1_p(Q)}.
\end{align}
Combining \eqref{eq3} and	\eqref{eq4},
  interpolating the $L^p$-norm
 of $Du$ and
choosing $\varepsilon>0$ to be	small enough
 we obtain
\begin{equation}\label{eq8}
 \|u\|_{W^{2,1}_p(Q)}\leq C(\|f\|_{\V_p(Q)}+\|\varphi\|_{\W_p(Q)}+
 \|u\|_{p,Q}).
\end{equation}

To estimate the  norm of $u$   we take into account once again	the
differential
equation (see \cite[Ch.~VII]{L}).  Let $\tl u$ be an extension of $u$ as a zero for
$(x,t)\not\in Q.$ Obviously $\|u\|_{p,Q}=\|\tl u\|_{p,\R^{n+1}}$ and the same
holds also for the time derivative. Further it	is easy to see that for any
$\varsigma\in (0,T)$
\begin{align*}
\int_{\R^n} |\tl u(x,\varsigma)|^p
dx&=\int_0^\varsigma\int_{\R^n}\frac{d}{dt}|\tl u(x,t)|^pdxdt\leq
p\int_0^\varsigma \int_{\R^n} |\tl u(x,t)|^{p-1}|\tl u_t(x,t)|dt\\
&\leq p\left(\int_0^\varsigma\int_{\R^n}|\tl u_t(x,t)|^pdxdt \right)^{1/p}
\left(\int_0^\varsigma\int_{\R^n}|\tl u(x,t)|^{p}dxdt \right)^{(p-1)/p}.
\end{align*}
According to the  above considerations and \eqref{eq8} we can write for the
first integral
$$
\|u_t\|_{p,\Omega\times(0,\varsigma)}\leq C  \|D^2u\|_{p,\Omega\times(0,\varsigma)}
+\|f\|^p_{p,Q}
 \leq	C \left(\|u\|_{p,\Omega\times(0,\varsigma)} +\|f\|_{\V_p(Q)}
+\|\varphi\|_{\W_p(Q)}	\right).
$$
The function $U(\varsigma)=\int_{\R^n}|\tl u(x,\varsigma)|^pdx$ satisfies
\begin{align*}
U(\varsigma)&\leq C\left(\int_0^\varsigma U(t)dt \right)^{(p-1)/p}
\left(\left(\int_0^\varsigma U(t)dt\right)^{1/p} +
\|f\|_{\V_p(Q)}+\|\varphi\|_{\W_p(Q)}  \right)\\
& \leq C \int_0^\varsigma U(t)dt   +
C\|u\|_{p,Q}^{p-1}(\|f\|_{\V_p(Q)}+\|\varphi\|_{\W_p(Q)}).
\end{align*}
Hence the classical  Gronwall inequality gives a bound for $U(\varsigma),$ that is
$$
U(\varsigma)\leq C \|u\|_{p,Q}^{p-1}(\|f\|_{\V_p(Q)}+\|\varphi\|_{\W_p(Q)})
$$
and  from the definition of $U$ it  follows
$$
\|u\|_{p,Q}\leq C ( \|f\|_{\V_p(Q)}+\|\varphi\|_{\W_p(Q)}).
$$
The last one combining with \eqref{eq8} and the a priori estimate for the solution 
of the Cauchy problem give exactly \eqref{eq5}.
$\Box$

\section{Unique solvability}\label{s4}

\noindent
{\bf Proof of Theorem~2.}
The {\it uniqueness\/} follows trivially from \eqref{eq5}. To prove 
  {\it solvability\/}	 we consider  perturbed problem
 \begin{equation}\label{2}
\begin{cases}
\P u^\varepsilon\equiv u^\varepsilon_t -
a^{ij}(x,t)D_{ij}u^\varepsilon(x,t)=f(x,t)& {\text in }\  Q,\\
\I u^\varepsilon \equiv u^\varepsilon(x,0)=0 &	{\text on}\	\Omega,\\
\B_\varepsilon u^\varepsilon \equiv\ds \frac{\partial
u^\varepsilon}{\partial\l_\varepsilon} =\varphi (x,t) & {\text on}\  S
\end{cases}
\end{equation}
where $\l_\varepsilon(x,t)=\l(x,t)+\varepsilon\nu(x).$ Obviously, for
$\varepsilon>0,$  $\l_\varepsilon$ is nowhere tangential to
 $S$ and whence the above problem, being regular, has unique
strong solution  $u^\varepsilon\in W^{2,1}_p(Q),$
$p\in(1,\infty).$ Moreover, in view of Theorem~\ref{th2}, $u^\varepsilon$
satisfies the estimate
\begin{equation}\label{eq16}
\|u^\varepsilon\|_{W^{2,1}_p(Q)}\leq C
 \left(\|f\|_{p,Q}+\|\varphi\|_{p,S}^{(1-1/p)} \right)\leq
C\left(\|f\|_{\V_p(Q)} +\|\varphi\|_{\W_p(Q)} \right).
\end{equation}
The second estimate is more restrictive but its constant does not depend on $\varepsilon$
while the first one depends on it throught the norm of $\l_\varepsilon$ as it  could 
be seen from  the a priori estimate in \cite{Sf}.

By virtue of the compactness of the embedding $W^{2,1}_p(Q)\hookrightarrow V^1_p(Q)\hookrightarrow L^p(Q),$ and the weak compactness of  bounded sets
in $W^{2,1}_p(Q)$ there exists
 subsequence, which we reable as $\{u^\varepsilon\},$
 converging
weakly to a function $u\in W^{2,1}_p(Q)$ and
	$\|u_\varepsilon-u\|_{V^1_p(Q)}\to 0$ as $\varepsilon\to 0.$ Since
$$
 \int_Q fgdxdt=\int_Q (\P u^\varepsilon)gdxdt\longrightarrow\int_Q (\P
u)gdxdt\qquad  {\rm as}\quad \varepsilon\to 0,\  \forall\  g\in L^{p/(p-1)}(Q)
$$
we must have $\P u=f$ a.e. in $Q.$
Extending $\l_\varepsilon$  and $\varphi$ in $Q$  preserving their regularity,
 we get
$$
\|\B_\varepsilon u^\varepsilon -\B_\varepsilon u\|_{p,Q}\leq
 C\|D(u_\varepsilon-u)\|_{p,Q}\to 0\quad {as}\quad \varepsilon\to 0.
$$
Hence
\begin{align*}
&\|\varphi-\B u\|_{p,Q}=\|\B_\varepsilon u^\varepsilon -\B u\|_{p,Q}\leq \|B_\varepsilon u^\varepsilon-\B_\varepsilon u\|_{p,Q}\\
 &+\|\B_\varepsilon u -\B u\|_{p,Q}\leq C
 \|D(u^\varepsilon-u)\|_{p,Q}+\varepsilon \|Du\|_{p,Q}\to 0\quad{\rm as}\
\varepsilon\to 0.
\end{align*}

Let we stress now our considerations on the cylinder $\Sigma.$ The solution of the 
regular problem gains two drivatives in $x$ from $f$ and one from $\varphi,$ i.e. 
$u^\varepsilon\in V^3_p(\Sigma).$ It could be seen easy taking the representation 
formula for the solution of \eqref{2} (see \cite[Lemma~1]{Sf}) and the estimates  
for the heat potentials 
(see \cite{LSU}) as it was done for the solution of the Cauchy problem in  Section~\ref{s6}. 
Than it halds 
$$
\|u^\varepsilon\|_{W^{2,1}_p(\Sigma)}+\|D^3_xu^\varepsilon\|_{p,\Sigma} +\|D_xD_tu^\varepsilon\|_{p,\Sigma}\leq C\big(\|f\|_{V^1_p(\Sigma)}+
\|\varphi\|_{p,\Sigma}^{2-1/p}\big).
$$
Repeating the above arguments we get $\|u^\varepsilon - u\|_{V^2_p(\Sigma)}\to 0$ as $\varepsilon\to 0.$ For the extensions of $\varphi$ and $\l_\varepsilon $
in $\Sigma$ we get 
$$
\|\B_\varepsilon u^\varepsilon -\B_\varepsilon u\|_{V^1_p(\Sigma)}\leq
 C\|u_\varepsilon-u\|_{V^2_p(\Sigma)}\to 0\quad {as}\quad \varepsilon\to 0
$$
and $\|\varphi-\B u\|_{V^1_p(\Sigma)}\to 0.$ 
Therefore $\varphi=\B u$ on $S$ in  trace sense.\quad
$\Box$


\begin{thebibliography}{MPSS}
\normalsize
\bibitem[1]{ADN}
Agmon, S.; Douglis, A.;	Nirenberg, L. {Estimates near the boundary for solutions of
elliptic
partial differential equations satisfying general boundary conditions},
{ Commun. Pure Appl. Math.} {\bf 1959}, {\it 12}, 623--727; ibid. {\bf 1964}, {\it 17},
35-92.
\bibitem[2]{BC}
 Bramanti, M.;  Cerutti, .C. $W_{p}^{1,2}$   solvability
for the Cauchy -- Dirichlet problem for parabolic equations with VMO
 coefficients, Commun. Part. Diff. Eq. {\bf  1993}, {\it 18},
1735--1763.
\bibitem[3]{E}
Egorov, Y.V. {\em Linear Differential Equations of Principal
     Type}, Contemporary Soviet Mathematics, New York 1986.
\bibitem[4]{GT}
 Gilbarg, D.;  Trudinger, N.S. {\em Elliptic Partial
     Differential Equations of Second Order}, 2nd~Ed., Springer--Verlag,
     Berlin, 1983.
\bibitem[5]{H}
H\"ormander, L. {\em The Analysis of Linear Partial Differential Operators,}
Springer-Verlag, Berlin 1983.
  \bibitem[6]{LSU}
 Ladyzhenskaya, O.A.;  Solonnikov, V.A.;   Ural'tseva, N.N.
    {\em Linear and Quasilinear Equations of Parabolic Type}, Transl. Math.
    Monographs, Vol.~23,	Amer. Math. Soc., Providence, R.I. 1968.
\bibitem[7]{L} Lieberman, G., {\em Second Order Parabolic Differential
Equations,} World Scientific Publishing, Singapore, 1996.
 \bibitem[8]{MPS}
 Maugeri, A.,  Palagachev, D.K.,  Softova, L.G. {\em  Elliptic and
Parabolic Equations with Discontinuous Coefficients},  Wiley-VCH, Berlin,
2000.
\bibitem[9]{MPV}
Maugeri, A.,  Palagachev, D.K.; Vitanza, C.
  {Oblique derivative problem  for uniformly elliptic
   operators with $VMO$ coefficients and applications,}
   {C. R. Acad. Sci. Paris, S\'er. I, Math.}
    {\bf  1998}, {\it 327}, 53--58.
\bibitem[10]{Mz}
Maz'ya, V.,  On a degenerate problem with directional derivative, { Math. USSR Sbornik}
 {\bf 1972}, {\it 16}, 429-469.
\bibitem[11]{Pl}
Palagachev, D.K. The tangential oblique derivative problem for second order
 quasilinear parabolic operators, { Comm. Part. Differ. Eq.} {\bf 1992}, {\it 17}, 867--903.
\bibitem[12]{P}
Paneah, B.P. {\em The Oblique Derivative Problem. The Poincar\'e Problem,}
	Wiley--VCH, Berlin 2000.
\bibitem[13]{PP}  Popivanov, P.R.;  Palagachev, D.K. {\em The Degenerate
Oblique Derivative Problem for Elliptic and Parabolic Equations}, Wiley-VCH
(Akademie-Verlag), Berlin, 1997.
  \bibitem[14]{S}
  Sarason, D.	 Functions of vanishing mean oscillation, {	Trans. Amer.
Math. Soc.} {\bf 1975}, {\it 207}, 391--405.
\bibitem[15]{Sf}
 Softova, L.G.  Oblique derivative problem for parabolic operators with
$VMO$	   coefficients, { Manuscr. Math.}
{\bf  2000}, {\it 103}, 203--220.
\bibitem[16]{W1} Winzell, B.  The oblique derivative problem I, { Math.
Ann.} {\bf 1977}, {\it 229}, 267-278.
\bibitem[17]{W2}
Winzell, B.  The oblique derivative problem II, {  Ark.
Math.}	{\bf 1979}, {\it 17}, 107-122.
\bibitem[18]{W3}
Winzell, B.  A boundary value problem with an oblique derivative,
{ Commun. Part. Diff. Eq.}	{\bf 1981}, {\it 6}, 305-328.
\end{thebibliography}
\end{document}